\documentclass[11pt]{article}
 \usepackage{amsmath,amssymb,amscd}
\newcommand{\pf}[1]{\trivlist \item[\hskip\labelsep\it #1\ ]}
\newcommand{\varpf}[1]{\trivlist \item[\hskip\labelsep\sc #1:]}
\newcommand{\qedbox}{$\rlap{$\sqcap$}\sqcup$}
\newcommand{\qed}{\qquad \qedbox \endtrivlist}
\newcommand{\varqed}{\hfill \rule{0.6em}{0.6em} \endtrivlist}

\newtheorem{theorem}{Theorem}[section]

\newtheorem{df}[theorem]{Definition}

\begin{document}
\title{Noether}
\author{\it Volker Runde}
\date{}
\maketitle
The highest honour that can be bestowed upon a mathematician is not the Fields Medal --- it is becoming an adjective: Euclid has been immortalized in Euclidean geometry, Descartes' memory is preserved in Cartesian coordinates, and Newton lives
on in Newtonian mechanics. Emmy Noether's linguistic monument are the {\it Noetherian rings\/}\footnote{I won't attempt to explain what they are, but you will encounter them when you take your first course in abstract algebra.}. To my knowledge, 
Emmy Noether is still the only female mathematician ever to have received this honour --- and she definitely was the first. 
\par
Amalie ``Emmy'' Noether was born on March 23, 1882, in the city of Erlangen, in the German province of Bavaria, as the first child of Max Noether and his wife Ida, n\'ee Kaufmann. The math gene ran in her family: Her father was a math professor at the 
University of Erlangen, and her younger brother Fritz would later become a mathematician, too.
\par
Getting an education that deserved the name was not easy for a woman in those days. In order to be formally enrolled at a German university you needed (and still need) the {\it Abitur\/}, a particular type of high school diploma, and in those days,
there were no schools that would allow girls to graduate with the {\it Abitur\/}. Emmy attended a {\it H\"ohere T\"ochter--Schule\/} in Erlangen, a school that provided the daughters of the bourgeoisie with an education that was deemed suitable for
girls, i.e.\ with an emphasis on languages and the fine arts. Science and mathematics were not taught there in any depth. After graduation in 1897, Emmy continued to study French and English privately, and three years later, the passed the
Bavarian state exam that allowed her to teach French and English at girls' schools. 
\par
Instead of working as a language teacher, Emmy spent the years between 1900 and 1903 auditing lectures at the University of Erlangen in subjects as history, philology, and --- of course! --- mathematics. During the same period, she also
started preparing for the {\it Abitur\/} exam. In July 1903, Emmy obtained her {\it Abitur\/}, and at about the same time, women who had the {\it Abitur\/} were allowed to officially enroll at Bavarian universities. After a year at G\"ottingen,
she enrolled at Erlangen in 1904, where she obtained her doctorate in 1907. Although, at that time, it was already quite an accomplishment for a woman to obtain a doctorate at all, her thesis did little to indicate her future greatness as a 
mathematician. Later in life, Emmy herself referred to her thesis in terms such as ``Rechnerei''\footnote{mere computations}, ``Formelgestr\"upp''\footnote{shrub of formulas}, and even ``Mist''\footnote{manure}. 
\par
Having received her doctorate, Emmy continued working as a mathematical researcher in Erlangen --- without a position and, of course, without pay. Nevertheless, in the following years she built a reputation as mathematician, and in 1909,
she was the first woman invited to speak at the annual congress of the German mathematical association. In 1915, the mathematicians Felix Klein and David Hilbert invited her to join them at G\"ottingen, in the province of Prussia. To say that, in those 
days, G\"ottingen was the center of the mathematical world, is an understatement. From the nineteenth century to the early 1930s, G\"ottingen simply {\it was\/} the mathematical world. Needless to say, that Emmy went.
\par
The invitation to G\"ottingen was not a job offer. There were serious legal obstacles to her becoming a professor. In Germany, a doctorate is not a sufficient qualification to become a professor: You need to obtain the {\it Habilitation\/}, which
is some kind of a second, but more demanding doctorate. Without doubt, Emmy Noether was a strong enough mathematician to obtain it, but Prussian law at that time only admitted males as candidates for the {\it Habilitation\/}. When Emmy, supported
by David Hilbert, nevertheless filed for her {\it Habilitation\/}, the Prussian ministry of culture intervened and forbade it. Despite this throwback, Emmy started teaching at G\"ottingen in 1915, albeit wihtout pay and not under her own name: 
Officially, Hilbert was the instructor, and she only assisted him.
\par
In 1919, many things changed in Europe. Monarchy in Germany was overthrown. For the first time, women had the right to vote, and Emmy Noether finally could file for the {\it Habilitation\/} without legal problems. This didn't mean that suddenly everything
had become easy for her (or for other women at German universities). Her {\it Habilitation\/} had to be approved of by the university's senate. The very thought of a woman receiving the {\it Habilitation\/} stirred a heated debate in the senate: 
If we grant her the {\it Habilitation\/} --- her opponents argued ---, then she might one day become a professor, and if she becomes a professor, she might be elected into the university's senate, and the idea of a woman sitting in the senate is
so abhorrent that it must not be allowed to happen. I'm not kidding: Such arguments were seriously brought forward against Emmy Noether's {\it Habilitation\/}. According to legend, it was David Hilbert who ended the debate in the senate with the 
memorable phrase: ``Meine Herren! Der Senat ist keine Badeanstalt!''\footnote{Gentlemen! The senate is not a public swimming pool!} As the first woman at a German university, Emmy Noether was granted the {\it Habilitation\/} in mathematics and 
finally could lecture under her own name (still without pay).
\par
In 1922, the University of G\"ottingen granted Emmy Noether the title of extraordinary professor. The extraordinary thing about this professorship was that it gave Emmy the right to call herself a professor, but did not come with any salary.
Until 1921, Max Noether had supported his daughter with money. When he died, her financial situation became tight. In order to prevent her from becoming destitute, the university finally gave her a paid teaching assignment in 1923: For the first
time in her life, at age 41, Emmy Noether, one of the leading mathematical researchers in Germany if not in the world, had an income of her own.
\par
The next ten years, were years of scholarly and professional success for her. At that time, a revolution took place in mathematics. Emphasis shifted from computations and explicit constructions to more abstract and conceptual approaches. 
David Hilbert was at the forefront of this revolution --- and so was Emmy Noether. In the 1920 and 1930, the field of algebra changed almost beyond recognition. Around 1900, algebra was about solving algebraic equations. Fifty years later,
it had become the study of algebraic structures such as groups, rings, and fields: It had become abstract algebra. One of the driving forces behind this shift was Emmy Noether. A former student of hers, Bartel van der Waerden, a Dutchman, later
wrote a textbook on abstract algebra which is, to a large extent, based on Emmy Noether's lectures at G\"ottingen. It has been a standard text for decades. 
\par
By all accounts, Emmy Noether was very popular with students. Although she was demanding and not a very good lecturer, she had a kind personality which compensated all that. She gathered a group of devoted followers --- nicknamed her satellites --- 
around her, with whom she went swimming in the municipal pool
and whom she invited to her small appartment for large bowls of pudding. Of course, besides swimming and eating pudding, Emmy and her satellites had long discussions about mathematics. It is interesting to see how her students dealt with her
being a woman: They chose to ignore it. For example, she was referred to as {\it der Noether\/}: In German, {\it der\/} is the definite article for the masculine gender.
\par
Every four years, mathematicians hold the International Congress of Mathematicians (ICM). To give an invited address at such a congress is a feat few mathematicians will ever accomplish. In 1932, Emmy Noether delivered an invited presentation at the
ICM in Zurich. As you may have guessed, she was the first woman to do so --- and for a long time after 1932 she would remain the only one. By the way, she still didn't have permanent employment at that time --- even though some of her former students
had already become professors.
\par
On January 30, 1933, Paul von Hindenburg, president of Germany, appointed Adolf Hitler as chancellor. A month later, the {\it Reichstag\/} went up in flames; the government put the blame on the communists. Three weeks later, the first concentration
camp began operating, and on April 1, the systematic persecution of the Jews started with a boycot of Jewish businesses and other facilities.
\par
Emmy Noether's political views were left wing. She espoused pacifism and, for some time, had been a member of the Social Democratic Party. And, more significantly, she was Jewish. On April 25, 1933, Emmy Noether was sent on leave, which, in fact,
meant that she had been fired.
\par
She chose not to take any risks and accepted a visiting position at Bryn Mawr women's college in Pennsylvania. The transition from working mainly as a researcher to undergraduate teaching must have been difficult for her, but she soon succeeded again
in surrounding herself with eager young minds. In 1934, Emmy returned to Germany only to cancel the lease of her appartment and to have her belongings shipped to America. A year later, she underwent brain surgery to have a tumor removed. Even close
friends hadn't known of her illness. On April 14, 1935, she died from complications following the operation. She still didn't have a permanent position.
\par
As I mentioned at the beginning of this article, Emmy's younger brother Fritz was also a mathematician. He being male, his career developed more smoothly. Eventually, he became a full professor in Breslau\footnote{now Wroc{\l}aw in Poland}. After the
Nazi's rise to power, he decided to leave Germany for the Soviet Union. There, he was one of the millions who disappeared in Stalin's reign of terror.
\par
By the way, if you want to see a Noetherian ring and don't want to wait till you take your first course in abstract algebra, there is an opportunity: There are Noetherian Rings at several univesities in North America, for instance at 
Berkeley\footnote{\tt http://www.math.berkeley.edu/$^\sim$nring/}, the University of Florida\footnote{\tt http://www.math.ufl.edu/$^\sim$nring/}, and the University of Wisconsin at 
Madison\footnote{\tt http://www.math.wisc.edu/$^\sim$hollings/noethring/}, all local organizations of women --- faculty and students, graduate and undergraduate, alike --- in mathematics.
\vfill
\begin{tabbing}
{\it Address\/}: \= Department of Mathematical and Statistical Sciences \\
\> University of Alberta \\
\> Edmonton, Alberta \\
\> Canada T6G 2G1 \\[\medskipamount]
{\it E-mail\/}: \> {\tt vrunde@ualberta.ca} \\[\medskipamount]
{\it URL\/}: \> {\tt http://www.math.ualberta.ca/$^\sim$runde/}
\end{tabbing} 
\end{document}